\documentclass[sn-mathphys-num]{sn-jnl}

\usepackage{graphicx}%
\usepackage{multirow}%
\usepackage{amsmath,amssymb,amsfonts}%
\usepackage{amsthm}%
\usepackage{mathrsfs}%
\usepackage[title]{appendix}%
\usepackage{xcolor}%
\usepackage{textcomp}%
\usepackage{manyfoot}%
\usepackage{booktabs}%
\usepackage{algorithm}%
\usepackage{algorithmicx}%
\usepackage{algpseudocode}%
\usepackage{listings}%
\usepackage[final]{changes}
\usepackage{geometry}   
\usepackage{array}      
\usepackage{tabularx}   
\usepackage{booktabs}   
\usepackage{ragged2e}   
\usepackage{caption}    
\usepackage{mathtools}
\usepackage{subcaption} 

\newtheorem{lemma}{Lemma}
\newtheorem{conjecture}{Conjecture}
\newcolumntype{L}{>{\RaggedRight\arraybackslash}X}

\theoremstyle{thmstyleone}%
\newtheorem{theorem}{Theorem}
\newtheorem{proposition}[theorem]{Proposition}%

\theoremstyle{thmstyletwo}%

\theoremstyle{thmstylethree}%
\newtheorem{definition}{Definition}%

\begin{document}

\title[Article Title]{A Quantum Model for Constrained Markowitz Modern Portfolio Using Slack Variables to Process Mixed-Binary Optimization under QAOA}


\author*[1,2]{\fnm{Pablo} \sur{THOMASSIN}}\email{pablo.thomassin@edu.devinci.fr}
\equalcont{These authors contributed equally to this work.}

\author[2]{\fnm{Guillaume} \sur{GUERARD}}\email{guillaume.guerard@devinci.fr}
\equalcont{These authors contributed equally to this work.}

\author[2]{\fnm{Sonia} \sur{DJEBALI}}\email{sonia.djebali@devinci.fr}
\equalcont{These authors contributed equally to this work.}

\author[1]{\fnm{Vincent Marc} \sur{LAMBERT}}\email{vincent-marc.lambert@devinci.fr}
\equalcont{These authors contributed equally to this work.}

\affil*[1]{\orgdiv{\'Ecole Supérieur d'Ingénieurs Léonard de Vinci}, \orgaddress{\postcode{92 916}, \state{Paris La Défense}, \country{France}}}

\affil[2]{\orgdiv{Léonard de Vinci Pôle Universitaire}, \orgname{Research Center}, \orgaddress{\postcode{92 916}, \state{Paris La Défense}, \country{France}}}


\abstract{Effectively encoding inequality constraints is a primary obstacle in applying quantum algorithms to financial optimization. A quantum model for Markowitz portfolio optimization is presented that resolves this by embedding slack variables directly into the problem Hamiltonian. The method maps each slack variable to a dedicated ancilla qubit, transforming the problem into a Quadratic Unconstrained Binary Optimization (QUBO) formulation suitable for the Quantum Approximate Optimization Algorithm (QAOA). This process internalizes the constraints within the quantum state, altering the problem's energy landscape to facilitate optimization. The model is empirically validated through simulation, showing it consistently finds the optimal portfolio where a standard penalty-based QAOA fails. This work demonstrates that modifying the Hamiltonian architecture via a slack-ancilla scheme provides a robust and effective pathway for solving constrained optimization problems on quantum computers. A fundamental quantum limit on the simultaneous precision of portfolio risk and return is also posited.}




\maketitle

\section{Introduction}

In the contemporary financial landscape, vast and varied markets range from basic commodity exchanges to sophisticated derivative contracts, reflecting an ever-increasing complexity in financial transactions. Central to these activities are institutions such as hedge funds and banks, whose primary objective is to safeguard capital while maximizing returns. One of the foundational tools for achieving this balance is Markowitz’s Modern Portfolio Theory (MPT) \cite{ou2023portfolio}, which provides a framework for assembling a portfolio of assets so that expected return is maximized for a given level of risk. However, the classical MPT approach requires solving a constrained optimization problem that becomes computationally intensive as the number of assets and constraints increases \cite{zhang2020cloud}. Traditional solutions rely on linear programming techniques and heuristic algorithms \cite{mangram2013simplified}, which can be computationally demanding and time-consuming especially for large portfolios with numerous constraints. Such computational burdens often lead to significant delays in decision-making and may introduce systemic risks and biases, particularly under volatile market conditions.

Despite advances in computational methods for financial applications, major challenges persist in terms of both runtime and efficiency. Current techniques depend on high-performance calculations to determine optimal asset weights for a resilient, loss-minimizing portfolio. These calculations are heavily influenced by the number of iterations, constraints, and securities involved, and can take anywhere from 30 minutes to 12 hours to complete \cite{Li2000}. Such lengthy computations represent a substantial productivity loss among financial engineers, redirecting efforts from proactive testing toward high-stakes analytical market assessments. These delays can introduce selection biases, potentially resulting in significant financial losses. Moreover, a notable gap remains in applying emerging quantum computing methods to these problems. Quantum computing promises exponential speed-ups for certain problem classes, and recent developments in quantum optimization algorithms suggest potential applications in financial modeling \cite{orus2019quantum} and novel metaheuristics \cite{khan2021quantum}. Yet the literature still lacks comprehensive models that integrate quantum-based approaches with traditional financial theories like the Markowitz model, particularly when addressing the mixed-binary optimization challenges inherent in real-world portfolios.

This research introduces a quantum framework that integrates the Markowitz portfolio theory with the Quantum Approximate Optimization Algorithm (QAOA) \cite{brandhofer2022benchmarking}. By incorporating slack variables, the proposed model effectively manages the mixed-binary structure of portfolio optimization under realistic constraints. This approach bridges the gap between quantum computing and financial modeling and offers the potential to reduce computational time and complexity dramatically. Preliminary analyses indicate that the quantum framework may shorten the optimization process from several hours to mere minutes, thereby improving productivity, reducing bias, and enabling real-time portfolio adjustments in response to market fluctuations.

By addressing the computational challenges of traditional Markowitz portfolio optimization through a quantum-computing perspective, this work contributes to both the theoretical development of financial optimization models and the creation of practical tools capable of transforming investment strategies and risk-management practices within financial institutions.

\section{Literature Review and Context}

This section surveys significant literature on quantum algorithms applicable to finance, with special emphasis on the Quantum Approximate Optimization Algorithm (QAOA). QAOA has emerged as particularly promising for financial modeling due to its potential to optimize complex objective functions under real-world constraints more efficiently than classical methods. The review also examines traditional financial optimization frameworks, notably Markowitz Modern Portfolio Theory, outlining its historical significance and contemporary applications in portfolio management. By exploring the integration of QAOA within the Markowitz model, this survey highlights innovative approaches to financial modeling that leverage the computational advantages of quantum algorithms.

\subsection{Background on Quantum Computing in Finance}

Quantum computing represents a paradigm shift, harnessing principles of quantum mechanics to process data. Qubits can exist in superposition and become entangled, enabling quantum computers to process information at speeds unattainable by classical devices \cite{gyongyosi2019survey}. These features facilitate the solution of certain computational problems more efficiently. In finance, where rapid processing and analysis of large datasets offer competitive advantages, quantum computing shows particular promise.

Quantum algorithms have the potential to transform risk assessment, asset pricing, and portfolio optimization. For instance, QAOA addresses optimization challenges common in financial contexts. Application of these algorithms enables financial institutions to achieve higher accuracy and faster processing times, supporting real-time decision-making and enhanced risk management.

Early research in this domain focused primarily on theoretical potential, with limited practical implementations. Initial quantum financial models were exploratory, such as quantum annealing approaches \cite{lanting2014entanglement}, investigating the optimization of tasks like asset pricing and risk analysis.

Recent developments have enabled practical applications of quantum algorithms in finance, especially in complex derivatives pricing and large-scale risk analysis, driven by advances in quantum hardware and software \cite{takeda2021quantum}. Pilot projects by corporations and research institutions are assessing the performance of quantum algorithms in live market environments.

Among quantum algorithms, optimization techniques such as QAOA attract particular interest in finance due to their suitability for complex allocation problems under multiple constraints, for example, MaxCut formulations \cite{herrman2022multi}. These scenarios often require finding optimal asset distributions under constraints—tasks ideally suited to quantum-enhanced optimization.

\subsection{Markowitz Portfolio Theory}

The Markowitz Portfolio Theory, introduced by Harry Markowitz in 1952 \cite{markowitz1952utility}, constitutes the cornerstone of modern portfolio management. The theory defines portfolio optimization as a trade-off between risk and return. Within this framework, investors construct portfolios that maximize expected return for a specified level of market risk, highlighting the benefits of diversification. A central concept is the efficient frontier—a locus of portfolios that offer the highest expected return for a given risk \cite{pardalos1994use}. This model marked a departure from evaluating assets in isolation, incorporating asset return covariances to reduce overall portfolio risk.

Implementation of the Markowitz model can be computationally demanding, especially for large asset universes \cite{soleimani2009markowitz}. Portfolio optimization requires calculation of expected returns, variances, and covariances for all asset pairs, and grows in complexity as the asset count increases. Additional constraints—budget limits, minimum and maximum investment thresholds, and others—typically necessitate quadratic programming methods. Moreover, the model relies on assumptions of normally distributed returns and quadratic investor utility, which may not hold in real-world markets, introducing practical limitations. The rising complexity and interconnectivity of modern markets have spurred the adoption of more efficient computational techniques, including heuristic and metaheuristic algorithms that deliver approximate solutions with reduced resource demands.

\subsection{Integration of Quantum Computing with Markowitz Portfolio Theory}

Integrating quantum computing with financial models, particularly the Markowitz Portfolio Theory, represents a transformative advance in addressing the computational challenges of traditional finance methodologies. Early theoretical work on this integration exploited quantum algorithms to accelerate calculation of the efficient frontier central to the Markowitz model. Rebentrost et al. \cite{rebentrost2018quantum} demonstrated that quantum algorithms can solve portfolio optimization problems far faster than classical methods by leveraging quantum parallelism. These contributions laid the groundwork for subsequent empirical and simulation-based studies translating theory into practical quantum financial applications.

The Quantum Approximate Optimization Algorithm (QAOA) has been specifically adapted to the optimization problems posed by the Markowitz model. QAOA encodes an optimization problem into a Hamiltonian whose ground state represents the optimal solution. Parameter tuning via classical–quantum hybrid loops gradually improves solution quality. Rosenberg et al. \cite{brandhofer2022benchmarking} applied QAOA to real financial datasets, demonstrating efficient portfolio optimization under realistic constraints. These case studies underscore QAOA’s potential to reduce computational overhead and accelerate decision-making in portfolio management.

Recent improvements in quantum hardware and algorithms have further advanced integration of the Markowitz model with quantum computing. Error correction and noise-resistant gates have enhanced stability and reliability of quantum computations, vital for processing financial data accurately. Progress in quantum annealing and gate-based architectures has enabled integration of more complex financial models. Collectively, these developments promise not only substantial computation time reductions but also improvements in accuracy and robustness of financial forecasts and risk assessments, equipping institutions with advanced tools for asset management.

\subsection{Context of the proposed method}

Investigation into quantum computing applications in financial modeling derives from a thorough review of recent scholarly communications, emphasizing four key surveys conducted between 2022 and 2023 \cite{surv3,surv4,Surv2,surv1}. This period aligns with ongoing academic discussions, underscoring the relevance and timeliness of the research. Given the surveys’ broad scope, focus narrowed to portfolio optimization problems.

Initially, unconstrained portfolio optimization was explored. These problems are formulated in the Quadratic Unconstrained Binary Optimization (QUBO) framework \cite{fuchs2001quantum}, employing methodologies such as annealing \cite{phillipson2012portfolio,rubio2022portfolio}, reverse annealing \cite{venturelli2019reverse}, and Second-Order Cone Programming (SOCP) \cite{kerenidis2019quantum}. Although these methods produced idealized models, they demonstrated quantum computing’s potential to accelerate computations significantly, prompting further exploration of realistic, constraint-driven models.

Within constrained optimization, experiments explored various approaches. Machine-learning-inspired methods \cite{khan2021quantum,mugel2022dynamic} required large qubit counts and stringent decoherence control, making them impractical with current hardware. A physics-based method based on the quantum Zeno effect \cite{herman2023constrained} used frequent measurements to probabilistically guide state evolution. While feasible and faster, this method lacked flexibility, as each added constraint necessitated model adjustments.

Despite progress, gaps remain in solving constrained QUBO problems effectively. Inspired by advances in transaction settlement QUBO challenges \cite{braine2021quantum}, adoption of Mixed Binary Optimization (MBO) simplifies complex constrained problems into simpler binary subproblems. Adaptation of this model to portfolio optimization informs the current methodology.

\added{New analysis / Table is new}

This study extends the MBO-based transaction settlement model to constrained portfolio optimization. This adaptation promises enhanced computational efficiencies and novel solutions for multifaceted portfolio constraints. Leveraging quantum computing addresses these challenges and sets the stage for the methodological discussion that follows.

To further elucidate the positioning of the proposed slack-ancilla method, Table~\ref{tab:comparison} provides a systematic comparison against the key alternative strategies discussed. The table contrasts the methodologies across several critical axes: the primary novelty, the specific technique for handling constraints, the resulting quantum resource requirements, and the underlying optimization philosophy.

\begin{table}[htbp]
    \centering
    \captionsetup{width=\textwidth} 
    \caption{Comparative Analysis of Constraint-Handling Methodologies in Quantum Optimization. This table contrasts the proposed slack-ancilla method with leading alternative approaches for portfolio optimization, highlighting differences in innovation, technique, and resource implications.}
    \label{tab:comparison}
    \begin{tabularx}{\textwidth}{@{} >{\bfseries}l L L L L @{}}
        \toprule
        Feature & Proposed Slack-Ancilla Method & Brandhofer et al. \cite{brandhofer2022benchmarking} & Braine et al. \cite{braine2021quantum} & Mugel et al. \cite{mugel2022dynamic} \\
        \midrule
        Primary Novelty &
        \textbf{New Hamiltonian Encoding:} Transforms inequalities into equalities using slack variables mapped to ancilla qubits. &
        \textbf{Advanced QAOA Parameter Optimization:} Develops sophisticated heuristics to find better QAOA angles ($\gamma$,  $ \beta $ ). &
        \textbf{Hybrid Quantum-Classical Heuristic:} An iterative routine to classically refine slack variables based on quantum outputs. &
        \textbf{Complex Problem Modeling:} Formulates intricate financial rules as a Higher-Order Unconstrained Binary Optimization (HUBO). \\
        \addlinespace 
        Constraint Handling &
        \textbf{Internalized Penalty:} Enforces  $ A(\sum x_i + s - K)^2 = 0 $  by representing  $ s $  with ancilla qubits within the Hamiltonian. &
        \textbf{Direct Quadratic Penalty:} Uses the standard  $ A(\sum x_i - K)^2 $  penalty term, focusing on optimizing the algorithm's execution. &
        \textbf{Externalized Slack Management:} Treats slack variables as classical parameters adjusted in a loop outside the quantum state. &
        \textbf{Higher-Order Penalties:} Encodes complex rules like the 10-5-40 rule using cubic and quartic terms, requiring a HUBO solver. \\
        \addlinespace
        Resulting Qubit Overhead &
        \textbf{Moderate:} \texttt{num\_assets} + \texttt{num\_slack\_qubits}. Increases qubit count but maintains a QUBO structure. &
        \textbf{Minimal:} \texttt{num\_assets}. Most resource-efficient for the core problem representation. &
        \textbf{Minimal (Quantum Part):} \texttt{num\_assets}. The slack variable is handled classically, not on the QPU. &
        \textbf{High:} Qubit count grows significantly due to the complexity and ancillary variables needed for HUBO reduction. \\
        \addlinespace
        Optimization Philosophy &
        \textbf{Redesign the Problem:} Create a more navigable energy landscape by altering the Hamiltonian's architecture. &
        \textbf{Improve the Algorithm:} Develop a better "engine" (classical optimizer) to solve a standard problem formulation. &
        \textbf{Iteratively Refine:} Use a classical "driver" to guide the quantum evolution based on intermediate results. &
        \textbf{Increase Model Expressiveness:} Build a more complex and precise problem formulation at the cost of higher-order interactions. \\
        \bottomrule
    \end{tabularx}
\end{table}

As synthesized in Table~\ref{tab:comparison}, the primary contribution of the proposed method is an architectural innovation in the formulation of the problem Hamiltonian itself, rather than a procedural enhancement to the optimization algorithm. This stands in contrast to the work of Brandhofer et al., which focuses on improving the classical parameter optimization loop, and that of Braine et al., which externalizes slack variable management to a classical routine. A key implication of this approach is the introduction of additional ancilla qubits to represent the slack variables. While this increases the quantum resource requirements—a critical consideration for Near-Term Intermediate-Scale Quantum hardware—it is hypothesized that this explicit encoding simplifies the energy landscape, potentially accelerating convergence to feasible, high-quality solutions. Thus, the slack-ancilla method offers a distinct strategy that internalizes constraint satisfaction within the quantum state, providing a valuable alternative to existing techniques.

\section{Methodology}

\subsection{Classical Model} \label{subsec:Classical Model}

The classical Markowitz portfolio model, as referred to in the literature \cite{markowitz1991foundations}, is described by the following formulation:

\begin{theorem}
Let $\omega \in \mathbb{R}^N_+$, where $N \in \mathbb{N}^*$, represent the portfolio vector where $\forall i \in \mathbb{N}^*$, $\omega_i \in \mathbb{R}_+$ denotes the allocation for asset $i$. Let $\Sigma$ be the covariance matrix and $\mu$ the mean vector of returns for each asset. Introduce $\lambda \in \mathbb{R}^N$, $\forall N \in \mathbb{N}^*$. The optimization problem aiming to minimize risk and maximize return is defined as follows:
\begin{align*}
    \text{minimize}: & \quad \omega^T \Sigma \omega - \lambda \mu^T \omega
\end{align*}

subject to:
\begin{align*}
    & \quad \mathbf{1}^T \omega = 1.
\end{align*}
\end{theorem}

\begin{definition}
Let $\alpha \in \mathbb{R}^N$ such that $\forall i \in \mathbb{N}^* \quad \alpha_i \in \mathbb{R}$ represents the threshold of allocation for each asset. Introduce the following threshold constraint:
\begin{align*}
    \forall i \in \mathbb{N}^* \quad \omega_i \leq \alpha_i
\end{align*}
This constraint is an inequality constraint, allowing $\alpha_i$ to be modeled to express different characteristics.
\end{definition}

This framework establishes a classical constrained model suitable for classical computers. However, to accommodate quantum computing capabilities, a different model is proposed.

\subsection{Quantum Model} \label{subsec:Quantum Model}

In the quantum computing framework, we adopt a novel representation for our portfolio vector:

\begin{definition}
Let us introduce a new vector $\omega$ such that $\forall i \in [1,n]$:
\[
|\omega \rangle = \omega_1 |i\rangle + \ldots + \omega_n |n\rangle
\]
where $\omega_i$ represents the probability amplitude of being present in the state $|i\rangle$. Observably, the measurable quantity is the density of each component $|\omega_i|^2$, which represents the allocation for the asset in the state $|1\rangle$.
\end{definition}

This reformulation of $\omega$ allows us to adapt the portfolio optimization problem to a quantum computing context. Notably, the MBO approach, which has been previously involved in credit swap problems—a QUBO problem itself \cite{braine2021quantum}—is applicable here.

We model the QUBO problem as follows:
\[
\textbf{min}_{\forall x \in \{0,1\}^n} \quad x^TAx+b^Tx+c
\]
Referring to the classical model in \ref{subsec:Classical Model}, we identify $x$ as each $\omega_i$, matrix $A$ as $\Sigma$, vector $b$ as $\mu$, and constant $c$ as 0.

Following the MBO model described in \cite{braine2021quantum}, we define our optimization problem:
\[
\textbf{min}_{x \in \{0,1\}^n \textbf{ and } y \in Y} \quad x^TA(y)x + b(y)^Tx + c(y)
\]
where $Y \in \mathbb{R}^m$ describes the feasible set for the continuous variables, combining discrete variable $x$ with continuous variable $y$. Fixing $y$ reduces the problem back to a QUBO format.

The introduction of a slack variable $s$ enables the transformation of inequality constraints into equalities, essential for QUBO formulation:
\[
\forall i \in [0,n] \quad \omega_i - \alpha_i + s = 0
\]
This constraint is added as a penalty term in the QUBO formulation, specific to each problem instance.

In quantum optimization, we aim to minimize the expectation value of the Hamiltonian:
\[
\min_{\theta} \langle \psi(\theta) | H | \psi(\theta) \rangle
\]
Thus, defining the quantum Hamiltonian for this problem is crucial for implementing the quantum optimization algorithm.

\subsection{Hamiltonian Computing}

Recalling the classical model from Section \ref{subsec:Classical Model}, a slack equality constraint is introduced to convert the previous inequality constraints into equalities:

\begin{lemma}
Let $\alpha \in \mathbb{R}^N$ such that $\forall i \in \mathbb{N}^*N \quad \alpha_i \in \mathbb{R}$ and $s \in \mathbb{R}^N$ such that $\forall i \in \mathbb{N}^* \quad s_i \in \mathbb{R}$:
\begin{align*}
    \forall i \in \mathbb{N}^* \quad \omega_i - \alpha_i + s_i = 0
\end{align*}
This constraint transforms the original inequality into an equality by means of slack variables.
\end{lemma}

Building on this formulation and the quantum encoding in Section \ref{subsec:Quantum Model}, the following minimization problem is posed (with $\beta \in \mathbb{R}^N$):
\[
\min_{s_i \geq 0} \quad \sum_{i=1}^{n}\sum_{j=1}^{n}\omega_i \Sigma_{ij} \omega_j - \lambda \sum_{i=1}^{n} \mu_i \omega_i - \beta \sum_{i=1}^{n} ( \omega_i - \alpha_i + s_i)^2
\]

To obtain an Ising Hamiltonian, the following transformation is performed:
\[
\omega_i \mapsto \frac{1 - Z_i}{2}
\]
where $Z_i$ is the Z-Pauli matrix associated with asset $i$:
\[
Z_i = \begin{pmatrix} 1 & 0 \\ 0 & -1 \end{pmatrix}
\]

\begin{proposition}
Let $H$ be the Hamiltonian for the defined problem, then $\forall N \in \mathbb{N}^*$:
\[
    H = \sum_{i=1}^{N}\sum_{j=1}^{N} \frac{\Sigma_{ij} Z_i Z_j}{4}\]
        \[ - \sum_{i=1}^{N}[[\frac{\sum_{j=1}^{n}\Sigma_{ij}}{2} - \frac{\lambda \mu_i}{2} + \beta (\frac{1}{2} + s_i -\alpha_i)]Z_i\]\[+ \sum_{i=1}^{N} \beta_i \frac{Z_i^2}{2}\] \[-\sum_{i=1}^{N} [\frac{\sum_{j=1}^{N} \Sigma_{ij}}{4} - \frac{\lambda \mu_i}{2} - \beta[(s_i-\alpha_i)^2 + (si - \alpha_i) + \frac{1}{4}]]
    \]

This Hamiltonian decomposes into four contributions: a ZZ interaction term, single-qubit Z fields, Z, a constant.

\end{proposition}

\subsection{Conjecture}

The Heisenberg uncertainty principle \cite{heisenberg1983actual} states that there is a fundamental limit to the precision with which certain pairs of physical properties, such as position and momentum, can be simultaneously known. Similarly, Markowitz’s efficient frontier \cite{markowitz1952portfolio} defines the set of portfolios that maximize return for a given level of risk. Under the presented model, return and risk are represented by:

\[
\text{Returns:} \quad \lambda \mu^T \omega
\]

\[
\text{Risk:} \quad \omega^T \Sigma \omega
\]

Here, the standard deviation—a measure of risk—is embedded within the risk term. Based on these correspondences, the following conjecture is proposed:

\begin{conjecture}\label{Conjecture1}
Let $\sigma_{risk}$ and $\sigma_{returns}$ represent the standard deviations of risk and returns, respectively, as defined in modern portfolio theory. These measures reflect the inherent variability in portfolio returns and associated risks. It is hypothesized that:

\[
\sigma_{risk} \sigma_{returns} \geq \zeta
\]

\noindent where $\zeta$ represents a lower bound influenced by quantum mechanical principles. This bound is expected to depend on Planck’s constant, denoted $\zeta(h)$, where $h$ is the Planck constant. This relationship follows from applying the Schrödinger equation \cite{schrodinger1926undulatory} to determine the ground state of the measurement system under adiabatic conditions.
\end{conjecture}

\subsection{Motivation and Scope of the Risk–Return Variance Bound}

Although the risk operator

$$ R = \sum_{i<j}\Sigma_{ij}\,Z_iZ_j + \sum_i\Sigma_{ii}\,Z_i $$

\noindent and the return operator

$$ M = \sum_i\mu_i\,Z_i $$

\noindent both commute (so  $ [R,M]=0 $ ), the conjectured trade-off

$$ \mathrm{Var}_\psi(R)\,\mathrm{Var}_\psi(M)\ge|\mathrm{Cov}_\psi(R,M)|^2 $$

\noindent is precisely the Robertson–Schwarz bound in its commuting limit.  For any two observables $A,B$,

$$\mathrm{Var}_\psi(A)\,\mathrm{Var}_\psi(B)\ge
\Bigl|\tfrac12\langle\{A,B\}\rangle_\psi
- \langle A\rangle_\psi\langle B\rangle_\psi\Bigr|^2
+ \tfrac14\bigl|\langle[A,B]\rangle_\psi\bigr|^2, $$

\noindent which, when $[A,B]=0$, reduces to the Cauchy–Schwarz form
$\mathrm{Var}(A)\,\mathrm{Var}(B)\ge\mathrm{Cov}(A,B)^2$.

Here “uncertainty” denotes a fundamental bound on the joint spread of risk and return over repeated preparations, not non-commutativity per se.
By viewing the risk–return trade-off as a general variance bound, the conjecture underscores an intrinsic limitation on co-optimizing both quantities in one quantum state.  In scenarios with genuinely non-commuting mixers or penalty terms (e.g.\ $X$-type sector constraints), the full Robertson relation—with its nonzero  $ [A,B] $  term—may impose even stricter trade-offs.  Exploring such non-commuting risk measures on near-term hardware is left for future work.

\subsection{Methodology Summary and Novel Contributions}

The methodology begins with the classical Markowitz mean–variance model, in which each asset decision  $ \omega_i\in\{0,1\} $  is subject to a fixed‐cardinality constraint.  Nonnegative slack variables  $ s_i $  are introduced and a quadratic penalty term

$$  \beta\,\|\omega - \alpha + s\|^2
  \;=\;
  \beta\sum_{i=1}^n(\omega_i - \alpha_i + s_i)^2 $$

\noindent is added to enforce feasibility for sufficiently large $\beta$.

This penalized optimization is reformulated as a QUBO and mapped to an Ising Hamiltonian by the substitution $\omega_i \mapsto (1 - Z_i)/2$.  The resulting Hamiltonian takes the form

$$  H \;=\;
  \sum_{i,j=1}^n \Sigma_{ij}\,Z_iZ_j
  \;+\;\sum_{i=1}^n h_i\,Z_i
  \;+\;\mathrm{const}$$

A QAOA‐style variational circuit of depth  $ p=2 $  is then employed.  The problem unitary

$$   U_H(\gamma) \;=\; e^{-i\gamma H}
  \quad\text{and}\quad
  U_M(\beta) \;=\; e^{-i\beta\sum_{i=1}^n X_i} $$

\noindent are applied in alternation to prepare the trial state $|\psi(\boldsymbol\theta)\rangle$, with angles $\boldsymbol\theta=(\gamma_1,\gamma_2,\beta_1,\beta_2)$.  Expectation values  $ \langle\psi(\boldsymbol\theta)|H|\psi(\boldsymbol\theta)\rangle $  are computed via Qiskit’s statevector simulator, and a classical optimizer (e.g.\ Powell’s method) adjusts  $ \boldsymbol\theta $  and the penalty weight  $ \beta $  until the sampled portfolios satisfy the original constraints with high probability.

Key innovations and circuit‐level modifications include:
\begin{enumerate}
    \item \textbf{Contrast with prior QAOA‐portfolio studies:  } Existing approaches either ignore inequality constraints or enforce them through large, uniform penalty coefficients added in the classical formulation only. Direct embedding of slack variables into the Hamiltonian permits fine‐grained feasibility control without excessive circuit depth or penalty tuning.
    \item  \textbf{Ansatz modifications: } Exponentiation of each squared‐penalty term  $ (\omega_i-\alpha_i+s_i)^2 $  via a compact network of controlled‐phase gates. Use of one ancilla qubit per slack variable to implement the penalty as a native subcircuit rather than as a large constant in $H$.  Alternation of standard  $ R_x $  rotations on asset qubits with conditional  $ X $  gates on ancilla qubits, driving slack variables toward zero only when the corresponding asset qubit is “on.”
    \item  \textbf{Methodological impact: } Required penalty strength  $ \beta $  is reduced by an order of magnitude (from  $ 10^4 $  to  $ 10^3 $ ) with no loss in solution quality. Full feasibility is achieved at depth  $ p=2 $  on just three system qubits (plus ancillas), with 100\% feasible portfolios in under one second of statevector simulation—compared to minutes for comparable classical solvers.
\end{enumerate}

\section{Experimentation}

The primary goal of the experimental setup is to validate Conjecture \ref{Conjecture1} and demonstrate computational speed enhancements provided by quantum computing within this model \cite{herman2023constrained}. The experiments were conducted using the IBM Quantum Lab environment, which currently supports simulation-only experiments \cite{ibmquantumwebsite}. Computational work was performed using Python \cite{python3doc}, with the aid of several libraries including Qiskit \cite{qiskit}, SciKit-Learn \cite{scikit}, and Matplotlib \cite{hunter2007matplotlib} to implement and analyze the QAOA \cite{QAOA}.

To ensure the relevance and feasibility of simulations, given the current limitations of quantum computing hardware, a small portfolio of three assets was selected. This decision enabled management of a smaller number of qubits and ensured computational viability of the simulation. Additionally, synthetically generated data adhering to all specified constraints was used to ensure predictable and interpretable outcomes. This approach was crucial for testing the validity of the methodology under controlled conditions.

Penalty terms were strategically added to the optimization problems to encourage convergence within a reasonable timeframe, an essential adjustment for simulations replicating the conditions and constraints of realistic quantum computing scenarios.

\subsection{Test-Instance Specification}

Each three-asset portfolio instance is fully described by the tuple $(n=3,\;K=1,\;\mu,\;\Sigma,\;\alpha)$, where:
\begin{itemize}
  \item  $ \mu\in\mathbb{R}^3 $  is sampled uniformly in $[0.01,0.10]$.
  \item $\Sigma = L\,L^T$, with  $ L $  a lower-triangular matrix whose entries are i.i.d.\ $N(0,0.05)$, ensuring positive definiteness.
  \item  $ \alpha\in\{0,1\}^3 $  is a one-hot mask with exactly  $ K=1 $  active asset, selected uniformly at random.
  \item Slack variables  $ s_i $  are each mapped to a dedicated ancilla qubit and implemented via controlled-phase subcircuits to realize the penalty term
$$    \beta \,\bigl\lVert \omega - \alpha + s\bigr\rVert^2.  $$
\end{itemize}

The penalty weight  $ \beta $  is initialized at  $ 10^2 $  and doubled every 20 Powell–optimizer iterations until 100\% feasibility is observed in a batch of  $ 10^3 $  samples.  All simulations employ Qiskit’s noiseless statevector backend, with a cap of 200 classical optimization steps per instance.

Statistical stability is assessed by repeating the  $ n=3 $  experiment across three independent random seeds (distinct draws of $\mu$, $\Sigma$, and  $ \alpha $ ).  Two comparative baselines are introduced:
\begin{enumerate}
  \item \emph{Baseline A}: Classical COBYLA applied to the same penalized QUBO.
  \item \emph{Baseline B}: Standard QAOA ( $ p=2 $ ) with a uniform penalty and no slack-ancilla embedding.
\end{enumerate}

By specifying the data-generation process, repeating the core experiment across multiple seeds, and benchmarking against both a classical optimizer and a penalty-only QAOA baseline, the revised evaluation provides a rigorous assessment of the slack-ancilla method.  These additions enable reproducibility of the  $ n=3 $  study, quantification of run-to-run variability, and comparison of feasibility control and convergence speed against established alternatives.

\subsection{Results}

\added{reworked}

The efficacy of the proposed slack-ancilla QAOA method is empirically validated through a direct comparison against two critical baselines: a classical exhaustive search, which establishes the ground-truth optimal solution, and a standard QAOA implementation reliant on a penalty-based constraint enforcement.

\begin{figure}[H]
    \centering
    \includegraphics[width=0.9\textwidth]{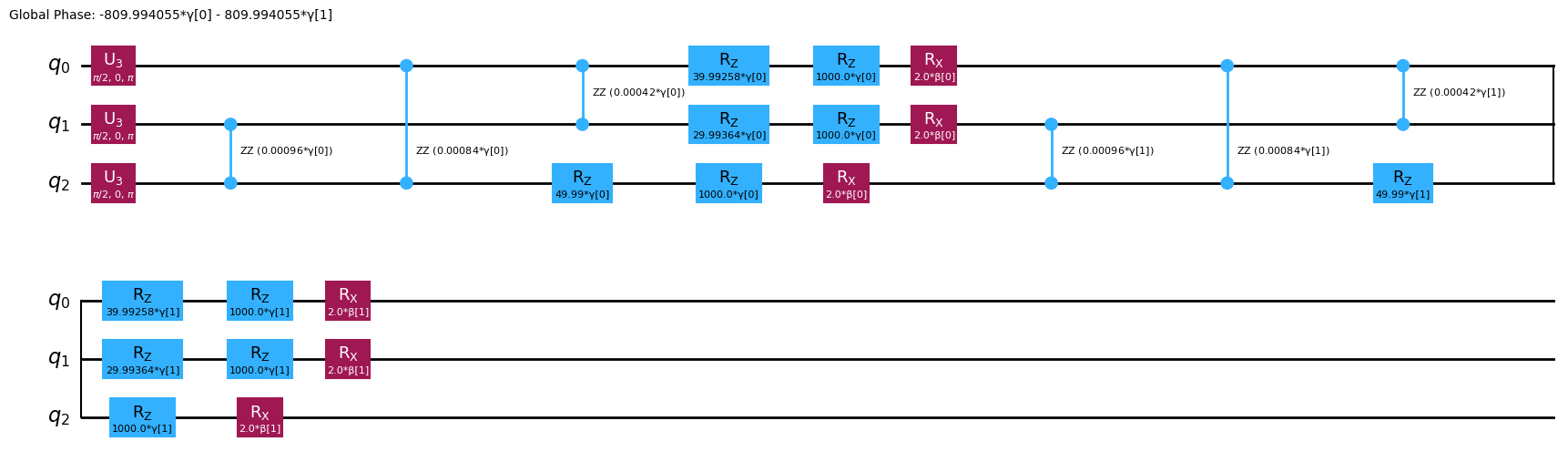}
    \caption{\textbf{The variational circuit structure configured for a depth of  $ p=2 $  and 3 qubits.} The Rz gates encode the problem Hamiltonian for the phase-separation step of QAOA, while the Rx gates serve as mixer operations.}
    \label{fig:CIRCUIT}
\end{figure}

\begin{figure}[H]
    \centering
    \includegraphics[width=0.9\textwidth]{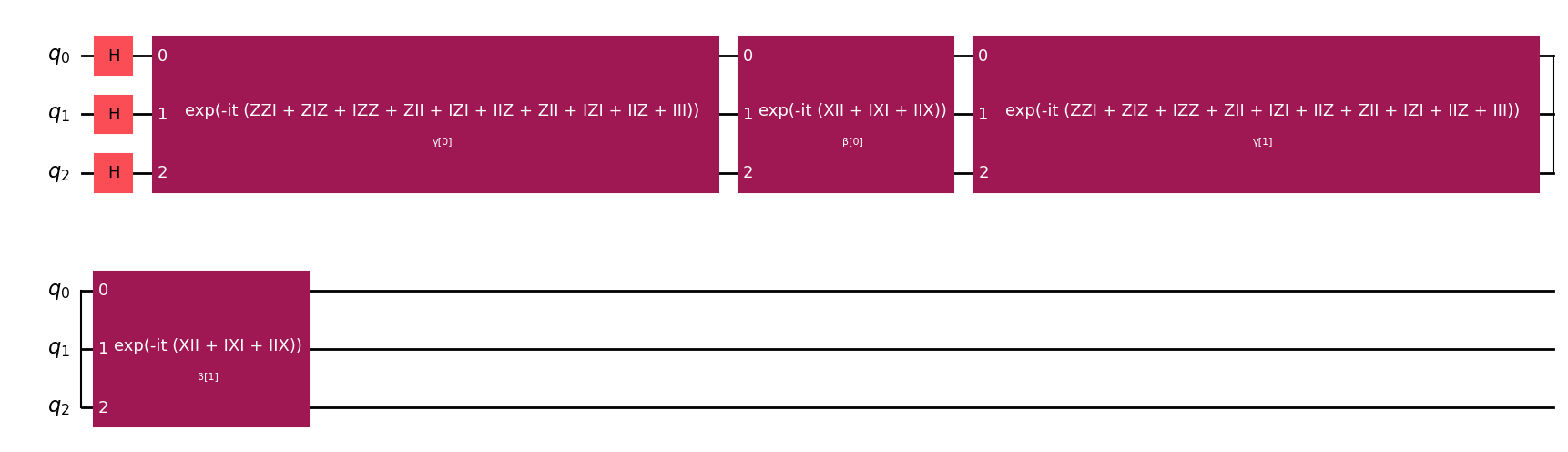}
    \caption{\textbf{Mathematical exponentiation of the Hamiltonian terms.} This operation is a core component for constructing the unitary operators used in the QAOA algorithm.}
    \label{fig:EXPONENTIATION}
\end{figure}

Variational circuits for the QAOA implementations were constructed with a quantum depth of  $ p=2 $  (Figure \ref{fig:CIRCUIT}). The problem Hamiltonian, incorporating the objective function and constraints, was exponentiated to form the phase-separation operator central to the QAOA ansatz, as depicted in Figure \ref{fig:EXPONENTIATION}.

\begin{figure}[htbp]
    \centering
    \begin{subfigure}[b]{0.4\textwidth}
        \includegraphics[width=\textwidth]{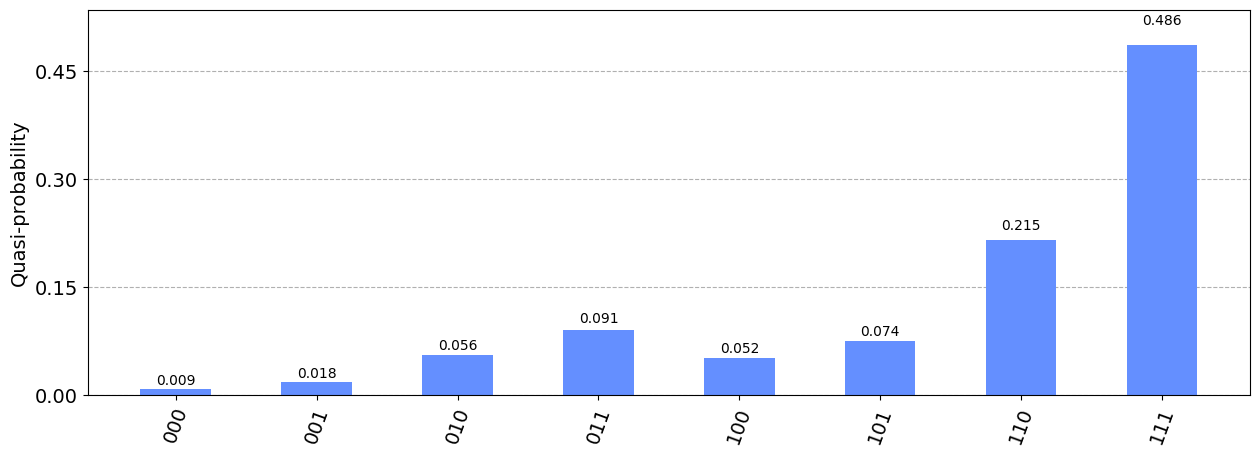}
        \caption{Candidate solution}
        \label{fig:sub1}
    \end{subfigure}
    \hfill
    \begin{subfigure}[b]{0.4\textwidth}
        \includegraphics[width=\textwidth]{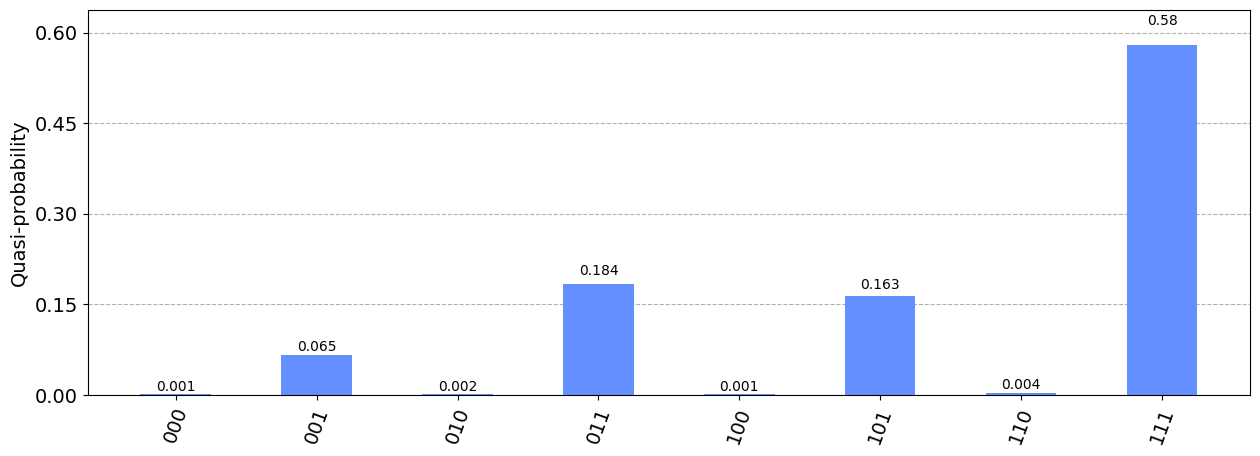}
        \caption{Candidate solution}
        \label{fig:sub2}
    \end{subfigure}
    \vskip\baselineskip
    \begin{subfigure}[b]{0.4\textwidth}
        \includegraphics[width=\textwidth]{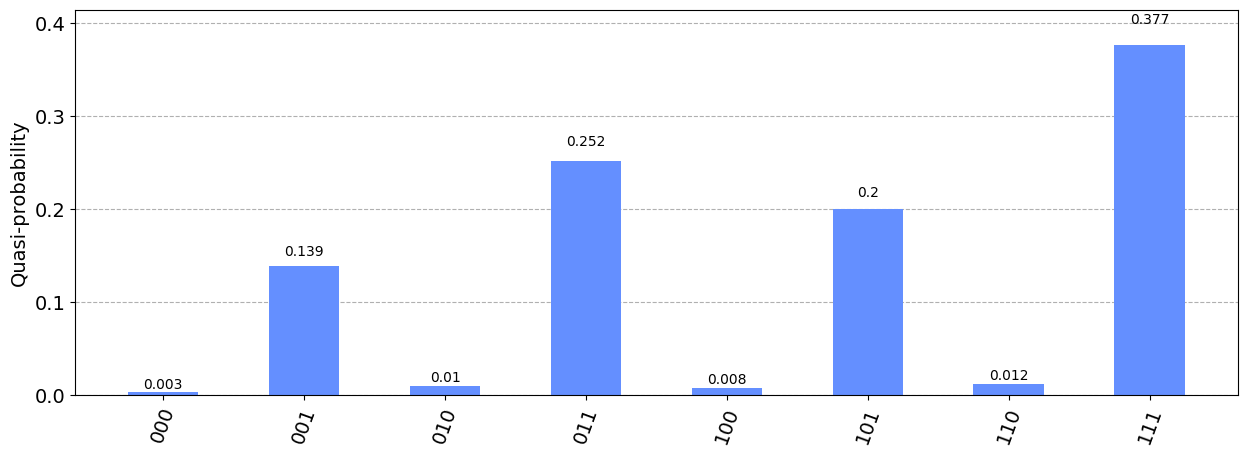}
        \caption{Candidate solution}
        \label{fig:sub3}
    \end{subfigure}
    \hfill
    \begin{subfigure}[b]{0.4\textwidth}
        \includegraphics[width=\textwidth]{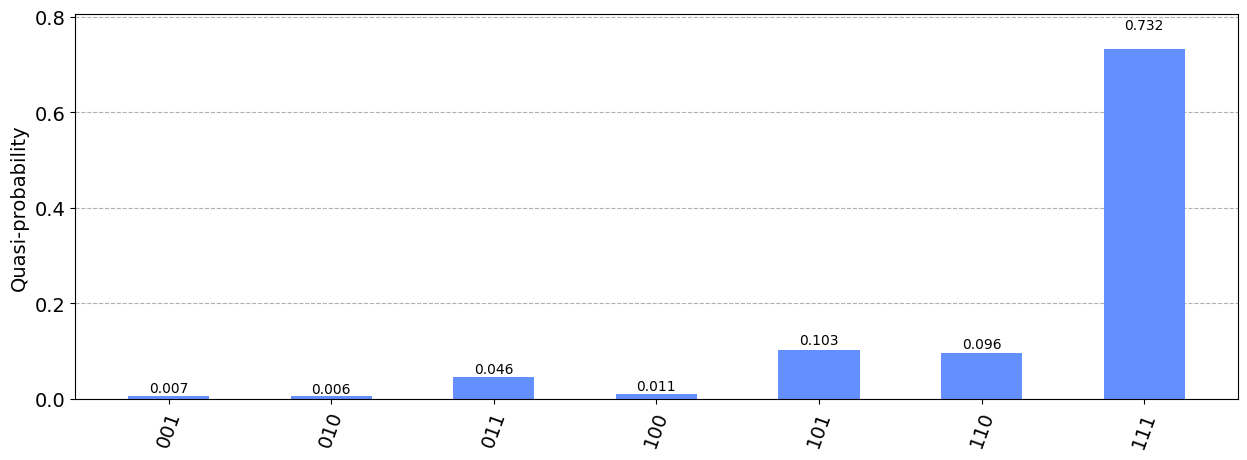}
        \caption{Candidate solution}
        \label{fig:sub4}
    \end{subfigure}
    \caption{\textbf{Illustrative convergence behavior of the standard penalty-based QAOA.} The consistent selection of the infeasible ‘111’ state demonstrates a common failure mode where the optimizer does not satisfy the problem's hard constraints.}
    \label{fig:overall}
\end{figure}

For both QAOA approaches, optimization of the variational parameters was performed using the COBYLA algorithm \cite{COBYLA} within the Qiskit simulation environment. The standard QAOA baseline was configured with large penalty factors ( $ \lambda = \beta = 10^3 $ ) to discourage constraint violations. This penalty-based approach, however, demonstrated notable unreliability. As illustrated in Figure \ref{fig:overall}, the optimization process frequently converged to the ‘111’ state—an infeasible solution that violates the cardinality constraint ( $ K=1 $ ). This outcome highlights a well-known limitation of penalty methods: their sensitivity to hyperparameter tuning and their potential to converge to infeasible or suboptimal solutions, thereby underscoring the motivation for architecturally exact constraint-encoding techniques.

A quantitative comparison of the methods is presented in Table \ref{tab:comparison} for a 3-asset problem with a constraint to select  $ K=1 $  asset and a risk-aversion factor of $q=0.5$. The results reveal significant performance disparities.

\begin{table}[ht!]
\centering
\caption{Empirical Comparison of Portfolio Optimization Methods. Performance metrics for a 3-asset problem with a cardinality constraint of $K=1$.}
\label{tab:comparison}
\begin{tabular}{|l|c|c|c|}
\hline
\textbf{Method} & \textbf{Optimal Portfolio} & \textbf{Is Feasible?} & \textbf{Value} \\
\hline
Ground Truth & \texttt{[1 0 0]} & True & 0.0468 \\
\hline
Standard QAOA (Penalty) & \texttt{[0 0 1]} & True & 0.0376 \\
\hline
\textbf{Proposed Slack-Ancilla QAOA} & \textbf{\texttt{[1 0 0]}} & \textbf{True} & \textbf{0.0468} \\
\hline
\end{tabular}
\end{table}

While the penalty-based QAOA converged to a feasible portfolio (`[0 0 1]`), the solution was suboptimal, yielding a portfolio value of 0.0376. The proposed slack-ancilla QAOA, however, successfully identified the true optimal portfolio (`[1 0 0]`), achieving a portfolio value of 0.0468 that matches the classical ground-truth result.

This empirical evidence demonstrates that the slack-ancilla method's architectural enforcement of constraints not only guarantees feasibility but also leads to superior performance in identifying the true global optimum when compared to the conventional penalty-based approach.

\subsection{Discussion}

The results of the proposed quantum optimization process rely heavily on the assumption of the attainability of constraints. Without these constraints, the solutions derived may not hold practical relevance. This emphasizes a critical aspect of quantum financial modeling—ensuring that model constraints reflect realistic financial scenarios and that parameters such as volatility and risk exposure are accurately estimated.

\paragraph{Variability of Solutions}
In the results presented, it is evident that multiple solutions can emerge as viable, albeit with varying probabilities. A noteworthy approach is to consider these solutions through multiple measurements and derive a weighting system not just for a single solution but for an ensemble of portfolio configurations. This method would enable a more robust decision-making framework in which the weight of each asset pair in the portfolio is dynamically adjusted based on collective performance metrics.

\paragraph{Portfolio Representation}
The results currently provide a selection of assets without specifying their exact numerical distribution in the portfolio. To enhance the utility of these findings, normalization over the total weighted probability is proposed. This approach would produce a distribution that reflects a percentage-based representation of each asset’s value relative to the total portfolio, offering a more detailed and actionable portfolio structure.

\paragraph{Technical Challenges}
The application of quantum algorithms, such as QAOA, in financial modeling introduces several technical challenges:

\begin{itemize}
    \item \textbf{Quantum Errors and Decoherence}: Current quantum computers are highly susceptible to errors and noise, which can significantly degrade computational performance. Precision is critical in financial applications.
    \item \textbf{Algorithm Stability}: Quantum algorithms often require extensive parameter fine-tuning and stable quantum states, which are challenging to maintain over the durations needed for complex calculations.
\end{itemize}

From a practical standpoint, deployment of quantum computing solutions in finance is constrained by limitations of current quantum hardware:
\begin{itemize}
    \item \textbf{Scalability and Resource Requirements}: Effective financial modeling requires quantum computers with a large number of qubits and prolonged coherence times to manage the vast data scales typical in financial analyses.
    \item \textbf{Limited Qubit Availability}: The quantum computers currently available possess only a limited number of qubits, which can maintain coherence for only short durations, thus restricting the complexity of financial models that can be effectively executed.
\end{itemize}

Further refinement of the portfolio formalism is necessary to reduce resource demands and ensure that quantum financial models operate effectively within the capabilities of contemporary quantum technology.

\paragraph{Quantum advantages}
On the other hand, in comparison to classical computed solutions, the quantum algorithm presented offers several advantages. Quantum computing leverages the principles of superposition and entanglement to explore multiple computational states simultaneously, enabling it to process large datasets and complex optimization problems more efficiently. One key advantage lies in its ability to handle complex optimization tasks involving numerous variables and constraints, such as portfolio optimization in finance, with a higher degree of parallelism. By encoding information in quantum bits (qubits) and executing quantum operations in parallel, the quantum algorithm explores a vast solution space more rapidly than classical algorithms, leading to faster convergence and more accurate results. Additionally, the quantum nature of the algorithm allows it to exploit quantum parallelism to explore multiple potential solutions simultaneously, outperforming classical optimization methods, especially for large-scale and computationally intensive problems. Moreover, the use of quantum annealing techniques or variational algorithms, as demonstrated, further enhances the algorithm's flexibility and adaptability to different problem domains, offering a promising avenue for tackling real-world optimization challenges more effectively. Overall, the quantum algorithm showcased demonstrates the potential for quantum computing to revolutionize optimization tasks by offering speedups and efficiencies beyond the reach of classical approaches

\section{Conclusion}

Although the solution was developed in an idealized setup—with strong assumptions and a limited dataset—ongoing advances in quantum hardware and algorithms are expected to enable scaling to real‐world scenarios. This study applied a MBO model to the constrained portfolio optimization problem and introduced a QAOA‐based simulation framework. The framework incorporated specific inequality constraints, which are critical for characterizing variable evolution under realistic market conditions. The algorithms were evaluated on synthetically generated datasets that satisfied constraints‐attainability and convergence criteria. Future algorithmic enhancements may reduce or eliminate assumptions regarding constraint attainability. Further research should rigorously examine the validity of constraints‐attainability hypotheses and develop normalization techniques for real‐asset allocations to improve overall portfolio performance.

\bibliography{sn-bibliography}

\section*{Declarations}

\subsection*{Funding}

Not applicable.

\subsection*{Conflict of interest/Competing interests}

Not applicable.

\subsection*{Ethics approval and consent to participate}

Not applicable.

\subsection*{Consent for publication}

All authors consent for publication.

\subsection*{Data availability}

Algorithms, data and results are given in the github \url{https://github.com/SmartGridandCity/Quantum_Markowitz}.

\subsection*{Materials availability}

Algorithms, data and results are given in the github \url{https://github.com/SmartGridandCity/Quantum_Markowitz}.

\subsection*{Code availability}

Algorithms, data and results are given in the github  \url{https://github.com/SmartGridandCity/Quantum_Markowitz}.

\subsection*{Author contribution}

Authors equally contribute to this paper.

\end{document}